\documentclass[10pt]{article}
\usepackage[utf8]{inputenc}
\usepackage{amsmath}
\usepackage{amssymb}
\usepackage{mathrsfs}
\usepackage{amsfonts}
\usepackage{mathrsfs,amscd,amssymb,amsthm,amsmath,bm,graphicx}
\usepackage{graphicx}
\usepackage{cite}
\usepackage{color}
\usepackage{enumitem}
\makeatletter

\renewcommand{\@seccntformat}[1]{{\csname the#1\endcsname}{\normalsize.}\hspace{.5em}}
\makeatother

\def \[{\begin{equation}}
\def \]{\end{equation}}

\newtheorem{thm}{Theorem}[section]

\newtheorem{lem}[thm]{Lemma}

\begin{document}
\setlength{\baselineskip}{13pt}
\begin{center}{\Large \bf The Normalized Laplacian Spectrum Analysis of Fractal M\"{o}bius Octagonal
Networks and its Applications
%\footnote{%The research is partially supported by National
%%Science Foundation of China (Grant No.  10671081)}
}

\vspace{4mm}

{\large Jia-Bao Liu $^{1}$}, ~Ting Zhang $^1$, ~Wenshui Lin $^{2,}$\footnote{Corresponding author. Email
address: wslin@xmu.edu.cn} \vspace{2mm}

{\small $^1$ School of Mathematics and Physics, Anhui Jianzhu University, Hefei 230601, P.R. China\\
\small $^2$ School of Informatics, Xiamen University, Xiamen 361005, China}
\vspace{2mm}
\end{center}

 {\noindent{\bf Abstract.}\ \
 The study and calculation of spectrum of networks can be used to describe networks structure and quantify
 analysis of networks performance.
The fractal M\"{o}bius octagonal networks, denoted by $Q_n$, is derived from the inverse identification of
the opposite lateral edges of fractal linear octagonal networks.
In this paper, the normalized Laplacian spectrum of $Q_n$
is determined by two  matrices $\mathcal {L}_A$ and $\mathcal {L}_S$. As an important application of our
results,  some topological indices (multiplicative degree-Kirchhoff index, the number of spanning trees)
formulas of $Q_n$ are obtained.

 \vspace{2mm}

\section{Introduction} \ \ \ \ \
In recent years, as a tool for studying mathematics, complex network plays a great role in scientific and
social research, which has attracted the attention of scholars and achieved some results \cite{D.J,A.L,WA}.
In this paper, the network is considered as simple, finite and undirected. If there is no special
explanation, the terminology and notation we mainly use come from \cite{J.A}.

A network can be regarded as a graph, $G=(V_{G},E_{G})$, and $V_G,~E_G$ are its vertex set and edge set,
respectively, $V_G=\{u_1,u_2,\ldots,u_n\},~E_G=\{e_1,e_2,\ldots,e_m\}$. The adjacent matrix $A_{G}$ of
simple networks can be expressed as
 \begin{eqnarray*}
(A_{G})_{ij}=\begin{cases}
         1, & if ~u_i~and~ u_j~ are~ adjacent; \\
         0, & otherwise.
\end{cases}
\end{eqnarray*}
The degree of $u_j$ is represented by $d_j$, for $1\leq j\leq n$. The Laplacian matrix $L_G$ of $G$ is
expressed as $L_G=D_G-A_G$,
where $D_G=diag_G(d_1,d_2,\ldots,d_n).$ Let $0=\mu_1<\mu_2\leq\cdots\leq \mu_n$ be the eigenvalues of $L_G$.
Then the Laplacian spectrum of $L_{G}$ is expressed as $Sp(L_{G})=\{\mu_1,\mu_2,\ldots,\mu_n\}.$

In the past few years, the normalized Laplacian matrix has come into the eyes of academic researchers. Some
of its results are not only useful for regular networks, but also suitable for general networks. For any
vertices $u_i~and~u_j$, the normalized Laplacian matrix can be written as
\begin{eqnarray*}
\big(\mathcal{L}(G)\big)_{ij}=
\begin{cases}
1,       & i=j; \\
-\frac{1}{\sqrt{d_id_j}},   & i\neq j,\ \ u_i\ \ and\ \ u_j\ \ are\ \ adjacent; \\
0,           &otherwise.
\end{cases}
\end{eqnarray*}

When each edge of a
connected network $G$ takes the place of a unit resistor\cite{b7}, the network can be regarded as a electric
circuit. $r_{ij}$ means the effective resistance distance between any vertices $u_i$ and $u_j$. It has been
found that this new parameter is also an invariant of network $G$, and it has also played an important role
in chemistry\cite{b8,b9}. At this time, based on the electrical network, Klein and Randi\'{c}\cite{b7}
proposed a new distance function called Kirchhoff index, $Kf(G)$. Similar to Wiener index, Kirchhoff index,
the sum of the effective resistance distances of any two vertices, is represented as
$Kf(G)=\sum_{i<j}r_{ij}$. Later, Gutman, Mohar\cite{b10} and Zhu et al.\cite{b11} introduced
\begin{eqnarray*}
Kf(G)=\sum_{i<j}r_{ij}=n\sum_{j=2}^{n}\frac{1}{\mu_j},
\end{eqnarray*}
where $0 = \mu_1 < \mu_2 < \cdots < \mu_n$ are the eigenvalues of $L_{G}$.

Chen and Zhang\cite{b34} proposed a new index, denoted by $DK(G),$ called multiplicative degree-Kirchhoff
index.
The mathematical expression is $DK(G)=\sum_{i<j}d_{i}d_{j}r_{ij}$. One can see\cite{b35,b36}. At the same
time, the index is closely related to the normalized Laplacian matrix of network $G$. For the studies of
normalized Laplacian of different networks, one can see \cite{b13,b14,b15,b16,b17}.

Furthermore, Kemeny's constant \cite{b18} is the expected number of time steps required for a Markov chain
to transition from a starting state $i$ to a random destination state sampled from the Markov chain's
stationary distribution, denoted as $Kc(G)$. The $Kc(G)$ is closely related to the effective resistance of
the graph, and the Kemeny's constant can be connected with the degree-Kirchhoff index by means of the
normalized Laplacian matrix.

Recently, a host of scholars have paid attention to study spectrum of networks. In 2016, Huang et al.
\cite{b25} characterized the normalized Laplacian spectrum, multiplicative degree-Kirchhoff index and the
number of spanning trees of linear hexagonal chain networks.
In 2018, Li et al.\cite{b19} studied the normalized Laplacian spectrum of the penta-graphene and M\"{o}bius
networks,
and obtained the corresponding formulas of degree-Kirchhoff index and the number of spanning trees. In 2018,
Li et al. \cite{b20} also studied the normalized Laplacian spectrum of linear phenylenes and their
dicyclobutadieno derivatives. By using the same methods, in 2019, X. Ma and H. Bian, \cite{b21,b22}
not only got the expressions for the M\"{o}bius networks, but also calculated some indices of cylinder
phenylene networks. More information about normalized spectrum. In 2019, Liu et al. \cite{b23} derived the
multiplication degree-Kirchhoff index and spanning tree of octagonal-quadrilateral, one can refer to
\cite{b24,b26,b27,b28}. Inspired by their
excellent works, we try to explore the expressions of multiplicative degree-Kirchhoff
index $DK(G)$ and the number of spanning trees $\tau(G)$ of fractal M\"{o}bius octagonal networks.

Let $L_n$ be a fractal linear octagonal networks of $n$ octagons, it is shown in Figure 1. Then, the fractal
M\"{o}bius octagonal networks $Q_n$ is obtained from $L_n$ by identifying the opposite lateral edges in
reversed way.

The rest of the work is arranged as follows. In Section 2, in order to get the results of this paper, we
introduce some theorems and terms.
In Section 3, closed-form formulae of degree-Kirchhoff index and the number of spanning trees are determined
for $Q_n$. In Section 4, we have done a full text summary.

\begin{figure}[htbp]
\centering\includegraphics[width=12cm,height=3cm]{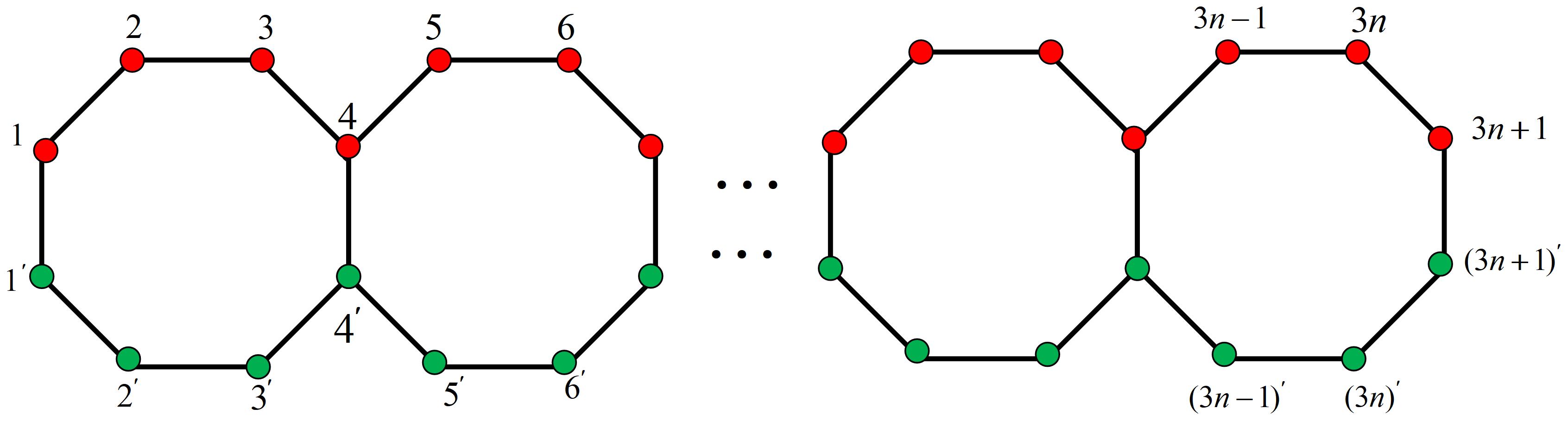}
\caption{ The fractal linear octagonal networks $L_n$.}\label{1}
\end{figure}

\begin{figure}[htbp]
\centering\includegraphics[width=12cm,height=4cm]{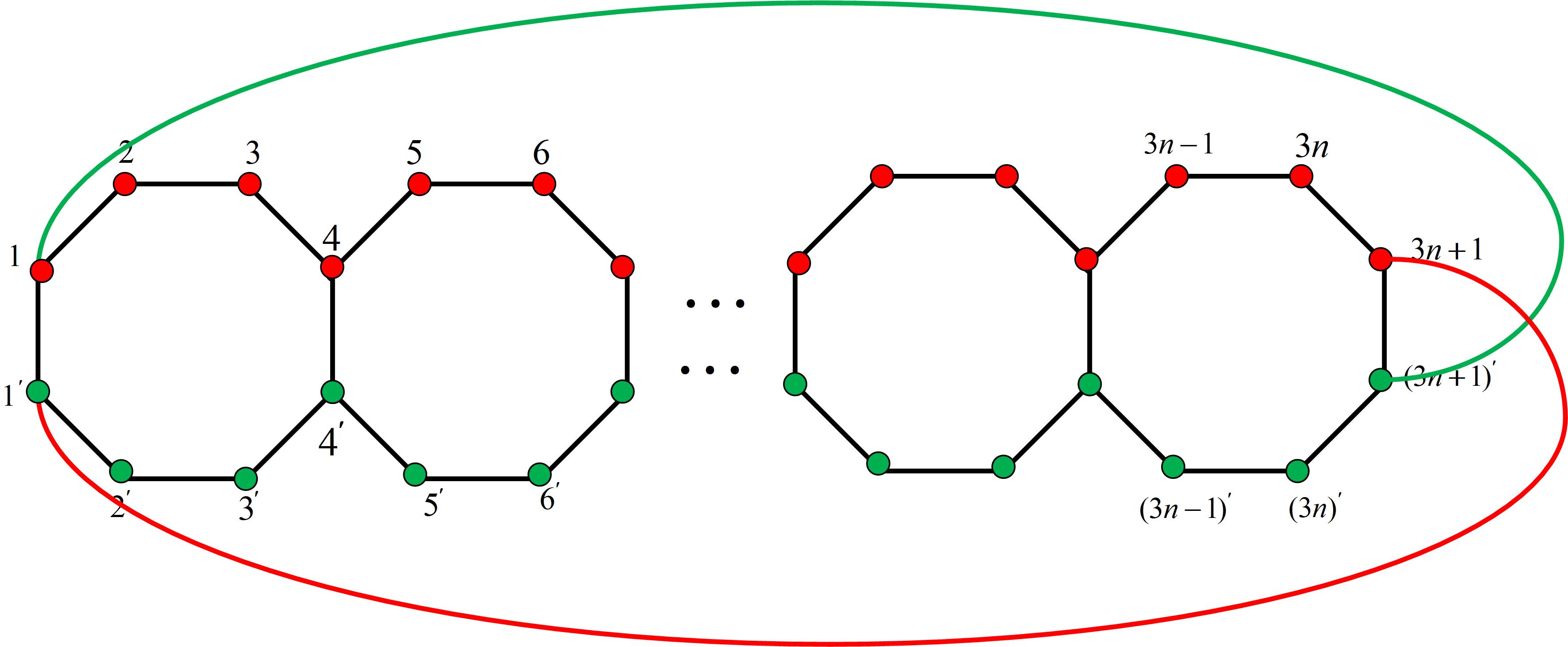}
\caption{ The fractal M\"{o}bius octagonal networks $Q_n$.}\label{2}
\end{figure}

\section{Preliminary}
\ \ \ \
In order to facilitate the calculation in the next section, firstly, we introduce some theorems and
lemmas. In this paper, we use $\Phi(G)=det(I-\mathcal{L}(G))$ as its characteristic polynomial of
$\mathcal{L}(G)$.
It is well known that the roots of $\Phi(G)$ are composed of the normalized Laplacian eigenvalues of
$\mathcal{L}(G)$.

Suppose $\pi=(\overline{1},\overline{2},\ldots,\overline{m})(1,1')(2,2')\cdots(n,n')$ is an automorphism of
a network $G=(V_{G},E_{G})$.
From Figure 2, obviously, we find that $|V_{G}|=m+2n$. Provided
$V_0=\{\overline{1},\overline{2},\ldots,\overline{m}\},V_1=\{1,2,\ldots,n\},V_2=\{1',2',\ldots,n'\}.$
Therefore,
$\mathcal L(G)$ can be written as
\begin{equation*}
\mathcal L(G)=\left(
  \begin{array}{ccc}
  \mathcal L_{V_0V_0}& \mathcal L_{V_0V_1} &\mathcal L_{V_0V_2}\\
\mathcal L_{V_1V_0}& \mathcal L_{V_1V_1} &\mathcal L_{V_1V_2}\\
\mathcal L_{V_2V_0}& \mathcal L_{V_2V_1} &\mathcal L_{V_2V_2}\\
  \end{array}
\right),
\end{equation*}
where $\mathcal L_{V_iV_j}$ are the matrix formed by the rows and columns of $\mathcal L(G)$ corresponding
to the vertices in $V_i \bigcup V_j, 0 \leq i, j \leq2.$

Suppose that
\begin{equation*}
U=\left(
  \begin{array}{ccc}
  I_{m} &0 &0\\
  0 &-\frac{1}{\sqrt{2}}I_{n}  &\frac{1}{\sqrt{2}}I_{n}\\
  0 &-\frac{1}{\sqrt{2}}I_{n} &\frac{1}{\sqrt{2}}I_{n}\\
  \end{array}
\right).
\end{equation*}

Consequently,
\begin{equation*}
UL(O_n)U'=\left(
  \begin{array}{ccc}
  \mathcal L_A& 0\\
    0& \mathcal L_{S}
  \end{array}
\right),
\end{equation*}
where $U'$ is the transposition of $U$, $\mathcal L_A=\left(
  \begin{array}{cc}
  \mathcal L_{V_0V_0}& \sqrt{2}\mathcal L_{V_0V_1}\\
    \sqrt{2}\mathcal L_{V_1V_0}& \mathcal L_{V_1V_1}+\mathcal L_{V_1V_2}
  \end{array}
\right),$ and $\mathcal L_{S}=\mathcal L_{V_1V_1}-\mathcal L_{V_1V_2}$.

Huang et al.\cite{b25} proposed the decomposition theorem of normalized Laplacian characteristic
polynomial.\\
\begin{thm} \cite{b25} Let $\mathcal L(G),~ \mathcal L_A,~and~ \mathcal L_S$ be defined as above.
Then we have
\begin{equation*}
\phi_{\mathcal L(G)}(z)=\phi_{\mathcal L_A}(z)\phi_{\mathcal L_S}(z).
\end{equation*}
\end{thm}

\begin{lem}\cite{b30} Let $G$ be connected graph with $|V_G|=n~and~|E_G|=m$. Then
$$DK(G)=2m\sum_{j=2}^{n}\frac{1}{\lambda_j},$$
where $0=\lambda_1<\lambda_2\leq \cdots \leq\lambda_n$ are eigenvalues of $\mathcal L(G)$.
\end{lem}

\begin{lem}\cite{b32} Let $G$ be $|V_G|=n~and~|E_G|=m$. Then
\begin{equation*}
Kc(G)=\sum_{j=2}^{n}\frac{1}{\lambda_j}.
\end{equation*}

\end{lem}
\begin{lem}\cite{b12} Let $\tau(G)$ be the number of spanning trees of a connected graph $G$ of order $n$
and m edges. Then
$$\prod_{j=1}^nd_j\prod_{j=2}^n\lambda_j=2m\tau(G).$$
\end{lem}

\section{Main results}
In this section, we mainly give some results and proofs. Using theorem 2.1 and Vieta's theorem, $DK(Q_n)$
and $\tau(Q_n)$ are determine.

Obviously, $|V(Q_n)|=6n,~and~\pi=(1,\overline{1})(2,\overline{2}) \cdots (3n,\overline{(3n)})$ is an
automorphism of $Q_n$. Let $V_0=\{\varnothing\},~V_1=\{1,2,...,3n\},$ and $V_2= \{\overline{1},
\overline{2},\ldots,\overline{(3n)}\}$. Thus $\mathcal L_A=\mathcal L_{V_1V_1}+\mathcal L_{V_1V_2}~ and~
\mathcal L_S=\mathcal L_{V_1V_1}-\mathcal L_{V_1V_2}.$ By elementary calculations, we have
\begin{eqnarray*}
 \mathcal {L}_{V_1V_1}&=&
\left(
  \begin{array}{cccccccccc}
1 & -\frac{1}{\sqrt{6}} & 0& 0& 0& \cdots& 0& 0& 0& 0\\
-\frac{1}{\sqrt{6}} & 1 & -\frac{1}{2}& 0& 0& \cdots& 0& 0& 0& 0\\
0& -\frac{1}{2} & 1 & -\frac{1}{\sqrt{6}} & 0& \cdots& 0& 0& 0& 0\\
0& 0& -\frac{1}{\sqrt{6}}&1 &-\frac{1}{\sqrt{6}} &  \cdots& 0& 0& 0& 0\\
0& 0& 0&-\frac{1}{\sqrt{6}} & 1 &\cdots  & 0& 0& 0& 0\\
\vdots& \vdots& \vdots& \vdots& \vdots&\ddots  & \vdots& \vdots& \vdots& \vdots\\
0& 0& 0& 0& 0& \cdots & 1&-\frac{1}{\sqrt{6}} & 0& 0\\
0& 0& 0& 0& 0& \cdots& -\frac{1}{\sqrt{6}}&1 &-\frac{1}{\sqrt{6}} & 0\\
0& 0& 0& 0& 0& \cdots& 0&-\frac{1}{\sqrt{6}} &1 &-\frac{1}{2}  \\
0& 0& 0& 0& 0& \cdots& 0& 0&-\frac{1}{2} &1  \\
  \end{array}
\right)_{(3n)\times (3n)},
\end{eqnarray*}
and
\begin{eqnarray*}
 \mathcal{L}_{V_1V_2}&=&
\left(
  \begin{array}{cccccccccc}
-\frac{1}{3} & 0 & 0& 0& 0& \cdots& 0& 0& 0& -\frac{1}{\sqrt{6}}\\
0 & 0 & 0& 0& 0& \cdots& 0& 0& 0& 0\\
0& 0 & 0 & 0 & 0& \cdots& 0& 0& 0& 0\\
0& 0&0 &-\frac{1}{3} &0 &  \cdots& 0& 0& 0& 0\\
0& 0& 0&0 & 0 &\cdots  & 0& 0& 0& 0\\
\vdots& \vdots& \vdots& \vdots& \vdots&\ddots  & \vdots& \vdots& \vdots& \vdots\\
0& 0& 0& 0& 0& \cdots & 0 &0 & 0& 0\\
0& 0& 0& 0& 0& \cdots& 0&-\frac{1}{3} &0 & 0\\
0& 0& 0& 0& 0& \cdots& 0&0 &0 &0  \\
-\frac{1}{\sqrt{6}}& 0& 0& 0& 0& \cdots& 0& 0&0 &0  \\
  \end{array}
\right)_{(3n)\times (3n)}.
\end{eqnarray*}
Hence
\begin{eqnarray*}
 \mathcal {L}_A&=&
\left(
  \begin{array}{cccccccccc}
\frac{2}{3} & -\frac{1}{\sqrt{6}} & 0& 0& 0& \cdots& 0& 0& 0& -\frac{1}{\sqrt{6}}\\
-\frac{1}{\sqrt{6}} & 1 & -\frac{1}{2}& 0& 0& \cdots& 0& 0& 0& 0\\
0& -\frac{1}{2} & 1 & -\frac{1}{\sqrt{6}} & 0& \cdots& 0& 0& 0& 0\\
0& 0& -\frac{1}{\sqrt{6}}&\frac{2}{3} &-\frac{1}{\sqrt{6}} &  \cdots& 0& 0& 0& 0\\
0& 0& 0&-\frac{1}{\sqrt{6}} & 1 &\cdots  & 0& 0& 0& 0\\
\vdots& \vdots& \vdots& \vdots& \vdots&\ddots  & \vdots& \vdots& \vdots& \vdots\\
0& 0& 0& 0& 0& \cdots & 1&-\frac{1}{\sqrt{6}} & 0& 0\\
0& 0& 0& 0& 0& \cdots& -\frac{1}{\sqrt{6}}& \frac{2}{3} &-\frac{1}{\sqrt{6}} & 0\\
0& 0& 0& 0& 0& \cdots& 0&-\frac{1}{\sqrt{6}} &1 &-\frac{1}{2}  \\
-\frac{1}{\sqrt{6}}& 0& 0& 0& 0& \cdots& 0& 0&-\frac{1}{2} &1  \\
  \end{array}
\right)_{(3n)\times (3n)},
\end{eqnarray*}
and
\begin{eqnarray*}
 \mathcal{L}_S&=&
\left(
  \begin{array}{cccccccccc}
\frac{4}{3} & -\frac{1}{\sqrt{6}} & 0& 0& 0& \cdots& 0& 0& 0& \frac{1}{\sqrt{6}}\\
-\frac{1}{\sqrt{6}} & 1 & \frac{1}{2}& 0& 0& \cdots& 0& 0& 0& 0\\
0& \frac{1}{2} & 1 & \frac{1}{\sqrt{6}} & 0& \cdots& 0& 0& 0& 0\\
0& 0& 0&-\frac{1}{\sqrt{6}} &\frac{4}{3} &  \cdots& 0& 0& 0& 0\\
0& 0& 0&0 & 1 &\cdots  & 0& 0& 0& 0\\
\vdots& \vdots& \vdots& \vdots& \vdots&\ddots  & \vdots& \vdots& \vdots& \vdots\\
0& 0& 0& 0& 0& \cdots & 1&-\frac{1}{\sqrt{6}} & 0& 0\\
0& 0& 0& 0& 0& \cdots& 0&\frac{4}{3} &-\frac{1}{\sqrt{6}} & 0\\
0& 0& 0& 0& 0& \cdots& 0&-\frac{1}{\sqrt{6}} &1 &-\frac{1}{2}  \\
\frac{1}{\sqrt{6}}& 0& 0& 0& 0& \cdots& 0& 0&-\frac{1}{2} &1  \\
  \end{array}
\right)_{(3n)\times (3n)}.
\end{eqnarray*}

Let $0=\alpha_1<\alpha_2\leq \cdots \leq\alpha_{3n}$ and $0<\rho_1\leq\rho_2\leq \ldots \leq\rho_{3n}$ be
the eigenvalues of $ \mathcal{L}_A$ and $\mathcal{L}_S,$
respectively. Hence $Sp(\mathcal L(Q_n))=\{\alpha_1,\alpha_2,..., \alpha_{3n},\rho_1,\rho_2, \ldots
,\rho_{3n}\}$.

\begin{thm} Let $0=\alpha_1 \leq \alpha_2 \leq \cdot\cdot\cdot \leq \alpha_{3n}$ be eigenvalues of $\mathcal
{L}_A$. Then
\begin{eqnarray*}
\sum_{j=2}^{3n}\frac{1}{\alpha_j}=\frac{147n^{2}-19}{84}.
\end{eqnarray*}
\end{thm}
\noindent\textbf{Proof.} Let $\phi_{\mathcal L_{A}}(z)=z^{3n}+d_1z^{3n-1}+ \cdots +d_{3n-1}z+d_{3n}$ be the
characteristic polynomial of $\mathcal L_A$. According to
Vieta's theorem, we have
\begin{eqnarray*}
\sum_{j=2}^{3n}\frac{1}{\alpha_j}=\frac{(-1)^{3n-2}d_{3n-2}}{(-1)^{3n-1}d_{3n-1}}.
\end{eqnarray*}
To determine $(-1)^{3n-2}d_{3n-2}$ and $(-1)^{3n-2}d_{3n-1}$, we need more preparations. Let
\begin{eqnarray*}
 \mathcal {L}_A^{0}&=&
\left(
  \begin{array}{cccccccccc}
\frac{2}{3} & -\frac{1}{\sqrt{6}} & 0& 0& 0& \cdots& 0& 0& 0& 0\\
-\frac{1}{\sqrt{6}} & 1 & -\frac{1}{2}& 0& 0& \cdots& 0& 0& 0& 0\\
0& -\frac{1}{2} & 1 & -\frac{1}{\sqrt{6}} & 0& \cdots& 0& 0& 0& 0\\
0& 0& -\frac{1}{\sqrt{6}}&\frac{2}{3} &-\frac{1}{\sqrt{6}} &  \cdots& 0& 0& 0& 0\\
0& 0& 0&-\frac{1}{\sqrt{6}} & 1 &\cdots  & 0& 0& 0& 0\\
\vdots& \vdots& \vdots& \vdots& \vdots&\ddots  & \vdots& \vdots& \vdots& \vdots\\
0& 0& 0& 0& 0& \cdots & 1&-\frac{1}{\sqrt{6}} & 0& 0\\
0& 0& 0& 0& 0& \cdots& -\frac{1}{\sqrt{6}}&\frac{2}{3} &-\frac{1}{\sqrt{6}} & 0\\
0& 0& 0& 0& 0& \cdots& 0&-\frac{1}{\sqrt{6}} &1 &-\frac{1}{2}  \\
0& 0& 0& 0& 0& \cdots& 0& 0&-\frac{1}{2} &1  \\
  \end{array}
\right)_{(3n)\times (3n)},
\end{eqnarray*}
\begin{eqnarray*}
 \mathcal {L}_A^{1}&=&
\left(
  \begin{array}{cccccccccc}
1 & -\frac{1}{2} & 0& 0& 0& \cdots& 0& 0& 0& 0\\
-\frac{1}{2} & 1 & -\frac{1}{\sqrt{6}}& 0& 0& \cdots& 0& 0& 0& 0\\
0& -\frac{1}{\sqrt{6}} & \frac{2}{3} & -\frac{1}{\sqrt{6}} & 0& \cdots& 0& 0& 0& 0\\
0& 0& -\frac{1}{\sqrt{6}}&1 &-\frac{1}{2} &  \cdots& 0& 0& 0& 0\\
0& 0& 0&-\frac{1}{2} & 1 &\cdots  & 0& 0& 0& 0\\
\vdots& \vdots& \vdots& \vdots& \vdots&\ddots  & \vdots& \vdots& \vdots& \vdots\\
0& 0& 0& 0& 0& \cdots & \frac{2}{3}&-\frac{1}{\sqrt{6}} & 0& 0\\
0& 0& 0& 0& 0& \cdots& -\frac{1}{\sqrt{6}}&1 &-\frac{1}{2} & 0\\
0& 0& 0& 0& 0& \cdots& 0&-\frac{1}{2} &1 &-\frac{1}{\sqrt{6}}  \\
0& 0& 0& 0& 0& \cdots& 0& 0&-\frac{1}{\sqrt{6}} &\frac{2}{3}  \\
  \end{array}
\right)_{(3n)\times (3n)},
\end{eqnarray*}
and
\begin{eqnarray*}
 \mathcal {L}_A^{2}&=&
\left(
  \begin{array}{cccccccccc}
1 & -\frac{1}{\sqrt{6}} & 0& 0& 0& \cdots& 0& 0& 0& 0\\
-\frac{1}{6} & \frac{2}{3} & -\frac{1}{\sqrt{6}}& 0& 0& \cdots& 0& 0& 0& 0\\
0& -\frac{1}{\sqrt{6}} & 1 & -\frac{1}{2} & 0& \cdots& 0& 0& 0& 0\\
0& 0& -\frac{1}{2}&1 &-\frac{1}{\sqrt{6}} &  \cdots& 0& 0& 0& 0\\
0& 0& 0&-\frac{1}{\sqrt{6}} & \frac{2}{3} &\cdots  & 0& 0& 0& 0\\
\vdots& \vdots& \vdots& \vdots& \vdots&\ddots  & \vdots& \vdots& \vdots& \vdots\\
0& 0& 0& 0& 0& \cdots & \frac{2}{3}&-\frac{1}{\sqrt{6}} & 0& 0\\
0& 0& 0& 0& 0& \cdots& -\frac{1}{2}&1 &-\frac{1}{\sqrt{6}} & 0\\
0& 0& 0& 0& 0& \cdots& 0&-\frac{1}{\sqrt{6}} &\frac{2}{3} &-\frac{1}{\sqrt{6}}  \\
0& 0& 0& 0& 0& \cdots& 0& 0&-\frac{1}{2} &1  \\
  \end{array}
\right)_{(3n)\times (3n)}.
\end{eqnarray*}
Let $W_j^{0}$ $(W_{j}^{1}/W_{j}^{2})$ be the sequential principal minor of order $i$ of $\mathcal
{L}_A^{0}$(resp. $\mathcal {L}_A^{1}/\mathcal {L}_A^{2}$). Take $w_j^{0}:=detW_j^{0}$, $w_{j}^{1}:=det
W_{j}^{1}$,$w_{j}^{2}:=det W_{j}^{2}$,
$w_0=w_{0}^{1}=w_{0}^{2}=1.$ Then we can get some results as below.\\

 \textbf{Fact 1.} For $1\leq j \leq 3n$,
 \begin{eqnarray*}
w_j^{0}&=&\begin{cases}
(1+j)\cdot(\frac{1}{12})^{\frac{j}{3}},  &\text{if}\ j \equiv0 \pmod 3;\\
(\frac{1}{3}+\frac{j}{3})\cdot(\frac{1}{12})^{\frac{j-1}{3}},  &\text{if}\ j \equiv1 \pmod 3;\\
(1+j)\cdot(\frac{1}{12})^{\frac{j-2}{3}}, &\text{if}\ j \equiv2 \pmod 3.
\end{cases}
\end{eqnarray*}

\textbf{Proof.} It's easy to get$$w_1^{0}=\frac{2}{3}, w_2^{0}=\frac{1}{2}, w_3^{0}=\frac{1}{3},
 w_4^{0}=\frac{5}{36}, w_5^{0}=\frac{1}{12}, w_6^{0}=\frac{7}{144},$$ and for $3\leq j\leq3n-1$,
\begin{eqnarray*}
w_j^{0}&=&\begin{cases}
w_{j-1}^{0}-\frac{1}{4}w_{j-2}^{0},  &\text{if}\ j \equiv0 \pmod 3;\\
\frac{2}{3}w_{j-1}^{0}-\frac{1}{6}w_{j-2}^{0},&\text{if}\ j \equiv1 \pmod 3;\\
w_{j-1}^{0}-\frac{1}{6}w_{j-2}^{0}, &\text{if}\ j \equiv2 \pmod 3.
\end{cases}
\end{eqnarray*}

For $1\leq j\leq n-1, let ~A_i=w_{3j}^{0},~ B_{j}=w_{3j+1}^{0},~ C_{j}=w_{3j+2}^{0}$. Then
$A_1=\frac{1}{3},~B_0=\frac{2}{3},~B_1=\frac{5}{36},
~C_0=\frac{1}{2},~C_1=\frac{1}{12}.$ For $j\geq 2$

\begin{eqnarray}
\begin{cases}
A_{j}=C_{j-1}-\frac{1}{4}B_{j-1}; \\
B_{j}=\frac{2}{3}A_{j}-\frac{1}{6}C_{j-1};\\
C_{j}=B_{j}-\frac{1}{6}A_{j}.
\end{cases}
\end{eqnarray}

 By substituting elimination method into $3.1,$ we have $A_j=2B_{j}+\frac{1}{12}B_{j-1}$, and
 $C_j=\frac{2}{3}B_{j}-\frac{1}{72}B_{j-1}.$
 Finally, we put the results into the second equation of $3.1$, and one has
 $$144B_j-24B_{j-1}+B_{j-2}=0.$$
 Thus
 $$B_j=(c_1+c_2j)(\frac{1}{12})^{j},$$
 $B_0$ and $B_1$ are introduced into the above formula.
\begin{eqnarray*}
\begin{cases}
c_{1}\cdot (\frac{1}{12})^{0}=\frac{2}{3};   \\
(c_1+c_2)\cdot (\frac{1}{12})=\frac{5}{36}.
\end{cases}
\end{eqnarray*}
So, $c_1=\frac{2}{3},~c_2=1.$

Thus,
\begin{eqnarray*}
\begin{cases}
A_{j}=(1+3j)\cdot (\frac{1}{12})^{j}; \\
B_{j}=(\frac{2}{3}+j)\cdot (\frac{1}{12})^{j};\\
C_{j}=(\frac{1}{2}+\frac{1}{2}j)\cdot (\frac{1}{12})^{j}.
\end{cases}
\end{eqnarray*}

The Fact 1 proved over.

In a similar way, we acquire Fact 2 and Fact 3.

 \textbf{Fact 2.} For $1\leq j \leq 3n$,
 \begin{eqnarray*}
w_j^{1}&=&\begin{cases}
(1+j)\cdot(\frac{1}{12})^{\frac{j}{3}},  &\text{if}\ j \equiv0 \pmod 3;\\
(\frac{1}{2}+\frac{j}{2})\cdot(\frac{1}{12})^{\frac{j-1}{3}},  &\text{if}\ j \equiv1 \pmod 3;\\
(\frac{1}{4}+\frac{j}{4})\cdot(\frac{1}{12})^{\frac{j-2}{3}}, &\text{if}\ j \equiv2 \pmod 3.
\end{cases}
\end{eqnarray*}

 \textbf{Fact 3.} For $1\leq j \leq 3n-1$,
 $$w_{j}^{2}=w_{j-1}^{0}-\frac{1}{6}w_{j-2}^{1}.$$

\textbf{Fact 4.} $(-1)^{3n-1}d_{3n-1}=21n^{2}(\frac{1}{12})^{n}.$\\

\textbf{Proof.}
\begin{eqnarray*}
(-1)^{3n-1}d_{3n-1}&=&\sum_{x=1}^{3n}det\mathcal {L}_A[x] \\
&=&\sum_{x=3,\ x\equiv0\ (\text{mod}\ 3)}^{3n}det\mathcal {L}_A[x]+\sum_{x=1,\ x\equiv1\ (\text{mod}\
3)}^{3n-2}det\mathcal {L}_A[x]+\sum_{x=2,\ x\equiv2\ (\text{mod}\ 3)}^{3n-1}det\mathcal {L}_A[x],
\end{eqnarray*}
where
\begin{eqnarray*}
det\mathcal {L}_A[x]&=&\begin{cases}
w_{x-1}^{0}\cdot w_{3n-x}^{0}-\frac{1}{6}w_{x-2}^{1}\cdot w_{3n-x-1}^{0},  &\text{if}\ x \equiv0 \pmod 3;\\
w_{x-1}^{0}\cdot w_{3n-x}^{1}-\frac{1}{6}w_{x-2}^{1}\cdot w_{3n-x-1}^{1},&\text{if}\ x \equiv1 \pmod 3;\\
w_{x-1}^{0}\cdot w_{3n-x}^{2}-\frac{1}{6}w_{x-2}^{1}\cdot w_{3n-x-1}^{2}, &\text{if}\ x \equiv2 \pmod 3.
\end{cases}
\end{eqnarray*}
By the Fact 1 and Fact 2, we obtain
\begin{eqnarray*}
\sum_{x=3,\ x\equiv0\ (\text{mod}\ 3)}^{3n}det\mathcal {L}_A[x]&=&w_{j-1}^{0}\cdot
w_{3n-x}^{0}-\frac{1}{6}w_{x-2}^{1}\cdot w_{3n-x-1}^{0}\\
&=&\sum_{j=3,\ x\equiv0\ (\text{mod}\ 3)}^{3n}[\frac{x}{6}(\frac{1}{12})^{\frac{x-3}{3}} \cdot
(1+3n-x)(\frac{1}{12})^{\frac{3n-x}{3}}-\\
&&\frac{1}{6} \cdot \frac{1}{2}(x-1)\cdot (\frac{1}{12})^{\frac{x-3}{3}} \cdot \frac{1}{6}(3n-x)\cdot
(\frac{1}{12})^{\frac{3n-x-3}{3}}]\\
&=&6n^{2}(\frac{1}{12})^{n}.
\end{eqnarray*}

In the same way, according to Fact $1-3,$ we get the following results.
\begin{eqnarray*}
\sum_{x=1,\ x\equiv1\ (\text{mod}\ 3)}^{3n-2}det\mathcal {L}_A[x]&=&9n^{2}(\frac{1}{12})^{n}.\\
\sum_{x=2,\ x\equiv2\ (\text{mod}\ 3)}^{3n-1}det\mathcal {L}_A[x]&=&6n^{2}(\frac{1}{12})^{n}.
\end{eqnarray*}
The desired result holds.\hfill\rule{1ex}{1ex}

%In what follows, we will concentrate on the computation of
%$(-1)^{3n-1}b_{3n-1}$.
\vskip 0.5 cm
\textbf{Fact 5.}
$(-1)^{3n-2}d_{3n-2}=\frac{147n^{4}-19n^{2}}{4}.$
\vskip 0.1 cm

\textbf{Proof.}
 $(-1)^{3n-2}d_{3n-2}$ is the sum of all principal minors obtained by deleting two rows and two columns of
 $\mathcal{L}_{A}$. So
\begin{eqnarray*}
(-1)^{3n-2}d_{3n-2}&=&\sum_{1\leq x\leq y}^{3n}det\mathcal{L}_{A}[x,y]\\
&=&M_0+M_1+M_2.
\end{eqnarray*}

Note that
\begin{eqnarray*}
M_0&=&\sum_{1\leq x< y\leq n}\det\mathcal{L}_{A}[3x,3y]+\sum_{1\leq x< y\leq
n-1}\det\mathcal{L}_{A}[3x,3y+1]+\sum_{1\leq x< y\leq n-1}\det\mathcal{L}_{A}[3x,3y+2];\\
M_1&=&\sum_{1\leq x< y\leq n}\det\mathcal{L}_{A}[3x+1,3y]+\sum_{1\leq x< y\leq
n-1}\det\mathcal{L}_{A}[3x+1,3y+1]+\sum_{1\leq x< y\leq n-1}\det\mathcal{L}_{A}[3x+1,3y+2];\\
M_2&=&\sum_{1\leq x< y\leq n}\det\mathcal{L}_{A}[3x+2,3y]+\sum_{1\leq x< y\leq
n-1}\det\mathcal{L}_{A}[3x+2,3y+1]+\sum_{1\leq x< y\leq n-1}\det\mathcal{L}_{A}[3x+2,3y+2].
\end{eqnarray*}

{\bf{Case 1.}} If $3x~(mod~3) \equiv0$, $3y~(mod~3) \equiv 0, 1 \leq x < y \leq n.$ That is,
\begin{eqnarray*}
 \sum_{1\leq x< y\leq n}\det\mathcal{L}_{A}[3x,3y]&=&\sum_{1\leq x< y\leq n}w_{x-1}^{0}\cdot
 w_{y-x-1}^{0}\cdot w_{3n-y}^{0}-\frac{1}{6} w_{x-2}^{1}\cdot w_{y-x-1}\cdot w_{3n-y-1}\\
 &=&\sum_{1\leq x< y\leq n}4(y-x)(3n-y+x)\cdot (\frac{1}{12})^{n}\\
 &=&3(n^{4}-n^{2})(\frac{1}{12})^{n}.
\end{eqnarray*}

{\bf{Case 2.}} For $1 \leq x < y \leq n-1,$
\begin{eqnarray*}
 \sum_{1\leq x< y\leq n-1}\det\mathcal{L}_{A}[3x,3y+1]&=&\sum_{1\leq x< y\leq n-1}w_{x-1}^{0}\cdot
 w_{y-x-1}^{0}\cdot w_{3n-y}^{1}-\frac{1}{6} w_{x-2}^{1}\cdot w_{y-x-1}^{0}\cdot w_{3n-y-1}^{1}\\
 &=&\sum_{1\leq x< y\leq n-1}6(y-x)(3n-y+x)\cdot (\frac{1}{12})^{n}\\
 &=&(\frac{9}{2}n^{4}-6n^{3}+\frac{3}{2}n^{2})(\frac{1}{12})^{n}.
\end{eqnarray*}

{\bf{Case 3.}} For $1 \leq x < y \leq n-1,$
\begin{eqnarray*}
 \sum_{1\leq x< y\leq n-1}\det\mathcal{L}_{A}[3x,3y+2]&=&\sum_{1\leq x< y\leq n-1}w_{x-1}^{0}\cdot
 w_{y-x-1}^{0}\cdot w_{3n-y}^{2}-\frac{1}{6} w_{x-2}^{1}\cdot w_{y-x-1}^{0}\cdot w_{3n-y-1}^{2}\\
 &=&\sum_{1\leq x< y\leq n-1}4(y-x)(3n-k+j)\cdot (\frac{1}{12})^{n}\\
 &=&(3n^{4}-2n^{3}+n^{2}-2n)(\frac{1}{12})^{n}.
\end{eqnarray*}

{\bf{Case 4.}} If $(3x+1)~(mod~3) \equiv1$, $3y~(mod~3) \equiv 0, 1 \leq x < y \leq n.$ So,
\begin{eqnarray*}
 \sum_{1\leq x< y\leq n}\det\mathcal{L}_{A}[3x+1,3y]&=&\sum_{1\leq x< y\leq n}w_{x-1}^{0}\cdot
 w_{y-x-1}^{1}\cdot w_{3n-y}^{0}-\frac{1}{6} w_{x-2}^{1}\cdot w_{y-x-1}^{1}\cdot w_{3n-y-1}^{0}\\
 &=&\sum_{1\leq x< y\leq n-1}6(y-x)(3n-y+x)\cdot (\frac{1}{12})^{n}\\
 &=&(\frac{9}{2}n^{4}+6n^{3}+\frac{3}{2}n^{2})(\frac{1}{12})^{n}.
\end{eqnarray*}

{\bf{Case 5.}} For $ 1 \leq x < y \leq n-1$,
\begin{eqnarray*}
 \sum_{1\leq x< y\leq n-1}\det\mathcal{L}_{A}[3x+1,3y+1]&=&\sum_{1\leq x< y\leq n-1}w_{x-1}^{0}\cdot
 w_{y-x-1}^{1}\cdot w_{3n-y}^{1}-\frac{1}{6} w_{x-2}^{1}\cdot w_{y-x-1}^{1}\cdot w_{3n-y-1}^{1}\\
 &=&\sum_{1\leq x< y\leq n-1}9(y-x)(3n-y+x)\cdot (\frac{1}{12})^{n}\\
 &=&\frac{27}{4}(n^{4}-n^{2})(\frac{1}{12})^{n}.
\end{eqnarray*}

{\bf{Case 6.}} For $1 \leq x < y \leq n-1,$
\begin{eqnarray*}
 \sum_{1\leq x< y\leq n-1}\det\mathcal{L}_{A}[3x+1,3y+2]&=&\sum_{1\leq x< y\leq n-1}w_{x-1}^{0}\cdot
 w_{y-x-1}^{1}\cdot w_{3n-y}^{2}-\frac{1}{6} w_{x-2}^{1}\cdot w_{y-x-1}^{1}\cdot w_{3n-y-1}^{2}\\
 &=&\sum_{1\leq x< y\leq n-1}6(y-x)(3n-y+x)\cdot (\frac{1}{12})^{n}\\
 &=&\frac{9}{2}(n^{4}+3n^{3}+\frac{3}{2}n^{2}+3n)(\frac{1}{12})^{n}.
\end{eqnarray*}

{\bf{Case 7.}} If $(3x+2)~(mod~3) \equiv2$, $3y~(mod~3) \equiv 0,$ for $ 1 \leq x < y \leq n,$
\begin{eqnarray*}
 \sum_{1\leq x< y\leq n}\det\mathcal{L}_{A}[3x+2,3y]&=&\sum_{1\leq x< y\leq n}w_{x-1}^{0}\cdot
 w_{y-x-1}^{2}\cdot w_{3n-y}^{2}-\frac{1}{6} w_{x-2}^{1}\cdot w_{y-x-1}^{2}\cdot w_{3n-y-1}^{0}\\
 &=&\sum_{1\leq x< y\leq n}4(y-x)(3n-y+x)\cdot (\frac{1}{12})^{n}\\
 &=&(3n^{4}+2n^{3}+n^{2}+2n)(\frac{1}{12})^{n}.
\end{eqnarray*}

{\bf{Case 8.}} For $1 \leq x < y \leq n-1,$
\begin{eqnarray*}
 \sum_{1\leq x< y\leq n-1}\det\mathcal{L}_{A}[3x+2,3y+1]&=&\sum_{1\leq x< y\leq n-1}w_{x-1}^{0}\cdot
 w_{y-x-1}^{2}\cdot w_{3n-y}^{1}-\frac{1}{6} w_{x-2}^{1}\cdot w_{y-x-1}^{2}\cdot w_{3n-y-1}^{1}\\
 &=&\sum_{1\leq x< y\leq n-1}6(y-x)(3n-y+x)\cdot (\frac{1}{12})^{n}\\
 &=&(\frac{9}{2}n^{4}-3n^{3}+\frac{3}{2}n^{2}-3n)(\frac{1}{12})^{n}.
\end{eqnarray*}

{\bf{Case 9.}} For $1 \leq x < y \leq n-1,$
\begin{eqnarray*}
 \sum_{1\leq x< y\leq n-1}\det\mathcal{L}_{A}[3x+2,3y+2]&=&\sum_{1\leq x< y\leq n-1}w_{x-1}^{0}\cdot
 w_{y-x-1}^{2}\cdot w_{3n-y}^{2}-\frac{1}{6} w_{x-2}^{1}\cdot w_{y-x-1}^{2}\cdot w_{3n-y-1}^{2}\\
 &=&\sum_{1\leq x< y\leq n-1}4(y-x)(3n-y+x)\cdot (\frac{1}{12})^{n}\\
 &=&3(n^{4}-n^{2})(\frac{1}{12})^{n}.
\end{eqnarray*}
The proof of Fact 5 completed.\hfill\rule{1ex}{1ex}

\begin{thm} Assume that $\rho_1< \rho_2\leq \cdots\leq \rho_{3n}$
are the eigenvalues of $\mathcal{L}_{S}$. One has
\begin{eqnarray}
\sum_{i=1}^{3n}\frac{1}{\rho_j}=\frac{\frac{37\sqrt{15}n}{30}\left[(\frac{4+\sqrt{15}}{12})^{n}-(\frac{4-\sqrt{15}}{12})^{n}\right]}
{(\frac{4+\sqrt{15}}{12})^{n}+(\frac{4-\sqrt{15}}{12})^{n}+2(\frac{1}{12})^{n}}.
\end{eqnarray}
\end{thm}
\noindent\textbf{Proof.} Let $\phi_{\mathcal L_{S}}(z)=z^{3n}+t_1z^{3n-1}+...+t_{3n-1}z+t_{3n}$ be the
characteristic polynomial of $\mathcal L_S$. So,
\begin{eqnarray*}
\sum_{j=1}^{3n}\frac{1}{\rho_j}=\frac{(-1)^{3n-1}t_{3n-1}}{(-1)^{3n-1}t_{3n}}=\frac{(-1)^{3n-1}t_{3n-1}}{det\mathcal
L_S}.
\end{eqnarray*}
To determine $(-1)^{3n-1}t_{3n-1}$, we need more preparations. Let
\begin{eqnarray*}
 \mathcal{L}_S^{0}&=&
\left(
  \begin{array}{cccccccccc}
\frac{4}{3} & -\frac{1}{\sqrt{6}} & 0& 0& 0& \cdots& 0& 0& 0& 0\\
-\frac{1}{\sqrt{6}} & 1 & \frac{1}{2}& 0& 0& \cdots& 0& 0& 0& 0\\
0& \frac{1}{2} & 1 & \frac{1}{\sqrt{6}} & 0& \cdots& 0& 0& 0& 0\\
0& 0& 0&-\frac{1}{\sqrt{6}} &\frac{4}{3} &  \cdots& 0& 0& 0& 0\\
0& 0& 0&0 & 1 &\cdots  & 0& 0& 0& 0\\
\vdots& \vdots& \vdots& \vdots& \vdots&\ddots  & \vdots& \vdots& \vdots& \vdots\\
0& 0& 0& 0& 0& \cdots & 1&-\frac{1}{\sqrt{6}} & 0& 0\\
0& 0& 0& 0& 0& \cdots& 0&\frac{4}{3} &-\frac{1}{\sqrt{6}} & 0\\
0& 0& 0& 0& 0& \cdots& 0&-\frac{1}{\sqrt{6}} &1 &-\frac{1}{2}  \\
0& 0& 0& 0& 0& \cdots& 0& 0&-\frac{1}{2} &1  \\
  \end{array}
\right)_{(3n)\times (3n)},
\end{eqnarray*}
and
\begin{eqnarray*}
 \mathcal{L}_S^{1}&=&
\left(
  \begin{array}{cccccccccc}
1 & -\frac{1}{2} & 0& 0& 0& \cdots& 0& 0& 0& 0\\
-\frac{1}{2} & 1 & -\frac{1}{\sqrt{6}}& 0& 0& \cdots& 0& 0& 0& 0\\
0& \frac{1}{\sqrt{6}} & \frac{4}{3} & -\frac{1}{\sqrt{6}} & 0& \cdots& 0& 0& 0& 0\\
0& 0& 0&-\frac{1}{\sqrt{6}} &1 &  \cdots& 0& 0& 0& 0\\
0& 0& 0&0 & 1 &\cdots  & 0& 0& 0& 0\\
\vdots& \vdots& \vdots& \vdots& \vdots&\ddots  & \vdots& \vdots& \vdots& \vdots\\
0& 0& 0& 0& 0& \cdots & 1&-\frac{1}{\sqrt{6}} & 0& 0\\
0& 0& 0& 0& 0& \cdots& 0&1 &-\frac{1}{2} & 0\\
0& 0& 0& 0& 0& \cdots& 0&-\frac{1}{2} &1 &-\frac{1}{\sqrt{6}}  \\
0& 0& 0& 0& 0& \cdots& 0& 0&-\frac{1}{\sqrt{6}} &\frac{4}{3}  \\
  \end{array}
\right)_{(3n)\times (3n)}.
\end{eqnarray*}
Let $Q_{i}^{0}$ $(Q_{i}^{1})$ be the sequential principal minor of order $i$ of $ \mathcal{L}_S^{0}$ (resp.
$\mathcal{L}_S^{1}$).
Take $q_{i}^{0}=detQ_{i}^{0},~q_{i}^{1}=detQ_{i}^{1}$. Then we  get the following results.

\textbf{Claim 1.} For $1\leq j \leq 3n$,
 \begin{eqnarray*}
q_j^{0}&=&\begin{cases}
(\frac{1}{2}+\frac{\sqrt{15}}{5})(\frac{4+\sqrt{15}}{12})^{\frac{j}{3}}+(\frac{1}{2}-\frac{\sqrt{15}}{5})(\frac{4-\sqrt{15}}{12})^{\frac{j}{3}},
&\text{if}\ j \equiv0 \pmod 3;\\
(\frac{2}{3}+\frac{17\sqrt{15}}{90})(\frac{4+\sqrt{15}}{12})^{\frac{j-1}{3}}+(\frac{2}{3}-\frac{17\sqrt{15}}{90})(\frac{4-\sqrt{15}}{12})^{\frac{j-1}{3}},
&\text{if}\ j \equiv1 \pmod 3;\\
(\frac{7}{12}+\frac{7\sqrt{15}}{45})(\frac{4+\sqrt{15}}{12})^{\frac{j-2}{3}}+(\frac{7}{12}-\frac{7\sqrt{15}}{45})(\frac{4-\sqrt{15}}{12})^{\frac{j-2}{3}},
&\text{if}\ j \equiv2 \pmod 3.
\end{cases}
\end{eqnarray*}

\textbf{Proof.} According to the specific calculation, one has
$q_1^{0}=\frac{4}{3},~q_2^{0}=\frac{7}{6},~q_3^{0}=\frac{5}{6},~q_4^{0}=\frac{11}{12},~q_5^{0}=\frac{7}{9}$.
For $1\leq j \leq 3n$,
\begin{eqnarray*}
q_j^{0}&=&\begin{cases}
q_{j-1}^{0}-\frac{1}{4}q_{i-2}^{0},  &\text{if}\ j \equiv0 \pmod 3;\\
\frac{4}{3}q_{j-1}^{0}-\frac{1}{6}q_{i-2}^{0},  &\text{if}\ j \equiv1 \pmod 3;\\
q_{j-1}^{0}-\frac{1}{6}q_{i-2}^{0}, &\text{if}\ j \equiv2 \pmod 3.
\end{cases}
\end{eqnarray*}

For $1\leq j \leq n$, Suppose that $e_j=q_{3j}^{0}.$ When $0\leq j \leq n-1,$ Suppose that
$f_{j}=q_{3j+1}^{0},$ and $g_{j}=q_{3j+2}^{0}.$
Then $e_1=\frac{5}{6},~f_0=\frac{4}{3},~f_1=\frac{11}{12},~g_1=\frac{7}{9}$. For $ j\geq 2$,
\begin{eqnarray}
\begin{cases}
e_{j}=g_{j-1}-\frac{1}{4}f_{j-1}; \\
f_{j}=\frac{4}{3}e_{j}-\frac{1}{6}g_{j-1};\\
g_{j}=f_{j}-\frac{1}{6}e_{j}.
\end{cases}
\end{eqnarray}

From the first and second expressions of 3.3, we can get $e_j=\frac{6}{7}f_j+\frac{1}{28}f_{j-1}.$
Then put the result into the third equation, we have $~g_j=\frac{6}{7}f_j-\frac{1}{168}f_{j-1},$ so
$g_{j-1}=\frac{6}{7}f_{j-1}-\frac{1}{168}f_{j-2}$.
Finally, substituting $e_j$ and $g_{j-1}$ into the second formula,
$$144f_j-96f_{j-1}+f_{j-2}=0.$$
So, $f_j=a_1(\frac{4+\sqrt{15}}{12})^{j}+a_2(\frac{4-\sqrt{15}}{12})^{j}$,
substituting the initial conditions $f_1$ and $f_2$ into the above formula,
$a_1=\frac{2}{3}+\frac{17\sqrt{15}}{90}, and ~a_2=\frac{2}{3}-\frac{17\sqrt{15}}{90}$ are obtained.
And then,
\begin{eqnarray*}
\begin{cases}
e_{j}=(\frac{1}{2}+\frac{\sqrt{15}}{5})(\frac{4+\sqrt{15}}{12})^{j}+(\frac{1}{2}-\frac{\sqrt{15}}{5})(\frac{4-\sqrt{15}}{12})^{j};
\\
f_{j}=(\frac{2}{3}+\frac{17\sqrt{15}}{90})(\frac{4+\sqrt{15}}{12})^{j}+(\frac{2}{3}-\frac{17\sqrt{15}}{90})(\frac{4-\sqrt{15}}{12})^{j};\\
g_{j}=(\frac{7}{12}+\frac{7\sqrt{15}}{45})(\frac{4+\sqrt{15}}{12})^{j}+(\frac{7}{12}-\frac{7\sqrt{15}}{45})(\frac{4-\sqrt{15}}{12})^{j},
\end{cases}
\end{eqnarray*}
as desired.

~~In the same way, we can get Claim 2. Here we omit the proof.

\textbf{Claim 2.} For $1\leq j \leq 3n$,
 \begin{eqnarray*}
q_j^{1}&=&\begin{cases}
(\frac{1}{2}+\frac{\sqrt{15}}{5})(\frac{4+\sqrt{15}}{12})^{\frac{j}{3}}+(\frac{1}{2}-\frac{\sqrt{15}}{5})(\frac{4-\sqrt{15}}{12})^{\frac{j}{3}},
&\text{if}\ j \equiv0 \pmod 3;\\
(\frac{1}{2}+\frac{3\sqrt{15}}{20})(\frac{4+\sqrt{15}}{12})^{\frac{j-1}{3}}+(\frac{1}{2}-\frac{3\sqrt{15}}{20})(\frac{4-\sqrt{15}}{12})^{\frac{j-1}{3}},
&\text{if}\ j \equiv1 \pmod 3;\\
(\frac{3}{8}+\frac{\sqrt{15}}{10})(\frac{4+\sqrt{15}}{12})^{\frac{j-2}{3}}+(\frac{3}{8}-\frac{\sqrt{15}}{10})(\frac{4-\sqrt{15}}{12})^{\frac{j-2}{3}},
&\text{if}\ j \equiv2 \pmod 3.
\end{cases}
\end{eqnarray*}

Using the properties of determinants, we have
\begin{eqnarray*}
\det \mathcal L_S&=& \left|
  \begin{array}{ccccccc}
    \frac{4}{3} & -\frac{1}{\sqrt{6}} & 0 & 0 & \cdots & 0 & \frac{1}{\sqrt{6}}\\
    -\frac{1}{\sqrt{6}} & 1 & -\frac{1}{2} & 0 & \cdots & 0 & 0\\
     0 & -1 & 2 &-1 & \cdots & 0 & 0\\
    0 & 0 & -1 & 4  & \cdots & 0 & 0\\
    \vdots & \vdots & \vdots & \vdots & \ddots & \vdots & \vdots\\
    0 & 0 & 0 & 0 & \cdots & 1 & -\frac{1}{2}\\
    \frac{1}{\sqrt{6}} & 0 & 0 & 0 & \cdots & -\frac{1}{2} & 1  \\
  \end{array}
\right|_{3n\times 3n}
 \end{eqnarray*}

 \begin{eqnarray*}
 ~~~~~~~~&=&\left|
  \begin{array}{ccccccc}
    \frac{4}{3} & -\frac{1}{\sqrt{6}} & 0 & 0 & \cdots & 0 & 0\\
    -\frac{1}{\sqrt{6}} & 1 & -\frac{1}{2} & 0 & \cdots & 0 & 0\\
     0 & -1 & 2 &-1 & \cdots & 0 & 0\\
    0 & 0 & -1 & 4  & \cdots & 0 & 0\\
    \vdots & \vdots & \vdots & \vdots & \ddots & \vdots & \vdots\\
    0 & 0 & 0 & 0 & \cdots & 1 & -\frac{1}{2}\\
    \frac{1}{\sqrt{6}} & 0 & 0 & 0 & \cdots & -\frac{1}{2} & 1  \\
  \end{array}
\right|_{3n\times 3n}
\end{eqnarray*}

 \begin{eqnarray*}
& &+\left|
  \begin{array}{ccccccc}
    \frac{4}{3} & -\frac{1}{\sqrt{6}} & 0 & 0 & \cdots & 0 & \frac{1}{\sqrt{6}}\\
    -\frac{1}{\sqrt{6}} & 1 & -\frac{1}{2} & 0 & \cdots & 0 & 0\\
     0 & -1 & 2 &-1 & \cdots & 0 & 0\\
    0 & 0 & -1 & 4  & \cdots & 0 & 0\\
    \vdots & \vdots & \vdots & \vdots & \ddots & \vdots & \vdots\\
    0 & 0 & 0 & 0 & \cdots & 1 & 0\\
    \frac{1}{\sqrt{6}} & 0 & 0 & 0 & \cdots & -\frac{1}{2} & 0  \\
  \end{array}
\right|_{3n\times 3n}\\
 &=&q_{3n}^{0}-\frac{1}{6}q_{3n-2}^{1}+2\cdot (\frac{1}{12})^{n}                \\
 &=&(\frac{4+\sqrt{15}}{12})^{n}+(\frac{4-\sqrt{15}}{12})^{n}+2(\frac{1}{12})^{n}.
    \end{eqnarray*}

Combined Claims 1 and 2, we have the following Claim.

\textbf{Claim 3.} $\det \mathcal
L_S=(\frac{4+\sqrt{15}}{12})^{n}+(\frac{4-\sqrt{15}}{12})^{n}+2(\frac{1}{12})^{n}.$

Next, we focus on calculating $(-1)^{3n-1}t_{3n-1}.$

\textbf{Claim 4.}
$(-1)^{3n-1}t_{3n-1}=\frac{37\sqrt{15}n}{30}\left[(\frac{4+\sqrt{15}}{12})^{n}-(\frac{4-\sqrt{15}}{12})^{n}\right].$

\textbf{Proof.} $(-1)^{3n-1}a_{3n-1}$ is the sum of all $3n-1$ order principal subexpression by deleting the
row by column of $\mathcal L_S$. One has
\begin{eqnarray}
(-1)^{3n-1}t_{3n-1}=\sum_{x=3,x\equiv0~(mod~3)}^{3n}det \mathcal
L_S[x]+\sum_{x=1,x\equiv1~(mod~3)}^{3n-2}det \mathcal L_S[j]+\sum_{x=2,x\equiv2~(mod~3)}^{3n-1}det \mathcal
L_S[x].
\end{eqnarray}

In the determinant of $det \mathcal L_S[x]$, we obtain

\textbf{Claim 5.}
 \begin{eqnarray*}
det\mathcal L_S[x]&=&\begin{cases}
\frac{7\sqrt{15}}{15}[(\frac{4+\sqrt{15}}{12})^{n}-(\frac{4-\sqrt{15}}{12})^{n}],  &\text{if}\ x \equiv0
\pmod 3;\\
\frac{3\sqrt{15}}{10}[(\frac{4+\sqrt{15}}{12})^{n}-(\frac{4-\sqrt{15}}{12})^{n}],  &\text{if}\ x \equiv1
\pmod 3;\\
\frac{7\sqrt{15}}{15}[(\frac{4+\sqrt{15}}{12})^{n}-(\frac{4-\sqrt{15}}{12})^{n}], &\text{if}\ x \equiv2
\pmod 3.
\end{cases}
\end{eqnarray*}

\textbf{Proof.} Let $det\mathcal L_S[x]=\left(
  \begin{array}{cc}
   L &M  \\
   N &O  \\
  \end{array}
\right)$, where $L$ is a square matrix of
order $(x-1)$, and $O$ is $(3n-x)\times(3n-x)$. Obviously
\begin{eqnarray} \left(
  \begin{array}{cc}
  0 &-I_{x-1}  \\
   I_{3n-x} &0  \\
  \end{array}
\right)^{T}     det\mathcal L_S[x]   \left(
  \begin{array}{cc}
  0 &-I_{x-1}  \\
   I_{3n-x} &0  \\
  \end{array}
\right) = \left(
  \begin{array}{cc}
   O &-N  \\
   -M &L  \\
  \end{array}
\right).
\end{eqnarray}

 Supposing that $\left(
  \begin{array}{cc}
  0 &-I_{x-1}  \\
   I_{3n-x} &0  \\
  \end{array}
\right)=S$. Then the above formula can be written as
\begin{eqnarray*}
S^{T} det\mathcal L_S[x] S= \left(
  \begin{array}{cc}
   O &-N  \\
   -M &L  \\
  \end{array}
\right).
\end{eqnarray*}

{\bf{Subcase 1.}} If $x\equiv0~(mod~3), ~3\leq x\leq 3n-3$,
\begin{eqnarray*}
S^{T} det\mathcal L_S[x] S = Q_{3n-1}^{0}=\mathcal L_S[3n]
\end{eqnarray*}
and
\begin{eqnarray*}\sum_{x=3,x\equiv0~(mod~3)}^{3n}det \mathcal L_S[x]&=&nq_{3n-1}^{0}\\
&=&n\left((\frac{7}{12}+\frac{7\sqrt{15}}{45})(\frac{4+\sqrt{15}}{12})^{n}+(\frac{7}{12}-\frac{7\sqrt{15}}{45})(\frac{4-\sqrt{15}}{12})^{n-1}\right)\\
&=&\frac{7\sqrt{15}n}{15}[(\frac{4+\sqrt{15}}{12})^{n}-(\frac{4-\sqrt{15}}{12})^{n}].
\end{eqnarray*}

Similarly,\\

{\bf{Subcase 2.}} If $x\equiv1~(mod~3), ~1\leq x\leq 3n-2$,
\begin{eqnarray*}
 S^{T} det\mathcal L_S[x] S = Q_{3n-1}^{1}=\mathcal L_S[1]
\end{eqnarray*}
and
\begin{eqnarray*}\sum_{x=1,x\equiv1~(mod~3)}^{3n-2}det \mathcal L_S[x]&=&nq_{3n-1}^{1}\\
&=&\frac{3\sqrt{15}n}{10}[(\frac{4+\sqrt{15}}{12})^{n}-(\frac{4-\sqrt{15}}{12})^{n}].
\end{eqnarray*}

{\bf{Subcase 3.}} If $x\equiv2~(mod~3), ~2\leq x\leq 3n-1$,
\begin{eqnarray}
  S^{T} det\mathcal L_S[x] S=\left(
  \begin{array}{cccccccccc}
1 & -\frac{1}{\sqrt{6}} & 0& 0& 0& \cdots& 0& 0& 0& 0\\
-\frac{1}{\sqrt{6}} & \frac{4}{3} & -\frac{1}{6}& 0& 0& \cdots& 0& 0& 0& 0\\
0& -\frac{1}{6} & 1 & -\frac{1}{2} & 0& \cdots& 0& 0& 0& 0\\
0& 0& 0&-\frac{1}{\sqrt{6}} &\frac{4}{3} &  \cdots& 0& 0& 0& 0\\
0& 0& 0&0 & -\frac{1}{\sqrt{6}} &\cdots  & 0& 0& 0& 0\\
\vdots& \vdots& \vdots& \vdots& \vdots&\ddots  & \vdots& \vdots& \vdots& \vdots\\
0& 0& 0& 0& 0& \cdots &\frac{4}{3}&-\frac{1}{\sqrt{6}} & 0& 0\\
0& 0& 0& 0& 0& \cdots& -\frac{1}{\sqrt{6}}&1 &-\frac{1}{2} & 0\\
0& 0& 0& 0& 0& \cdots& 0&-\frac{1}{\sqrt{6}} &1 &-\frac{1}{\sqrt{6}}  \\
0& 0& 0& 0& 0& \cdots& 0& 0&-\frac{1}{6} &\frac{4}{3} \\
  \end{array}
\right)_{(3n-1)\times (3n-1)}
\end{eqnarray}
and
\begin{eqnarray*}\sum_{x=2,x\equiv2~(mod~3)}^{3n-1}det \mathcal L_S[x]&=&nq_{3n-1}^{0}\\
&=&\frac{7\sqrt{15}n}{15}[(\frac{4+\sqrt{15}}{12})^{n}-(\frac{4-\sqrt{15}}{12})^{n}],
\end{eqnarray*}
as desired.

\begin{thm} Suppose $Q_n$ is a fractal M\"{o}bius octagonal networks of $n$ octagons. Then
\begin{eqnarray*}
Dk(Q_n)&=&14n(\sum_{j=2}^{3n}\frac{1}{\alpha_j}+\sum_{j=1}^{3n}\frac{1}{\rho_j})\\
&=&14n[(\frac{147n^{2}-19}{84})+\frac{\frac{37\sqrt{15}n}{30}[(\frac{4+\sqrt{15}}{12})^{n}-(\frac{4-\sqrt{15}}{12})^{n}]}{(\frac{4+\sqrt{15}}{12})^{n}+(\frac{4-\sqrt{15}}{12})^{n}+2(\frac{1}{12})^{n}}]\\
&=&\frac{147n^{3}}{6}-\frac{19n}{6}+14n\xi_n,
\end{eqnarray*}
where\\
\begin{eqnarray*}\xi_n&=&\frac{\frac{37\sqrt{15}n}{30}\left[(\frac{4+\sqrt{15}}{12})^{n}-(\frac{4-\sqrt{15}}{12})^{n}\right]}
{(\frac{4+\sqrt{15}}{12})^{n}+(\frac{4-\sqrt{15}}{12})^{n}+2(\frac{1}{12})^{n}}.
\end{eqnarray*}
\end{thm}

\begin{table}[h]
\setlength{\abovecaptionskip}{0.05cm} \centering\vspace{.3cm}
\caption{The Degree-Kirchhoff indices of $Q_n$ from $Q_{1}$
to $Q_{30}$.}
\begin{tabular}{cccccc}
  \hline
  $G$       &$DK(G)$       &$G$        &$DK(G)$        &$G$         &$DK(G)$ \\
  \hline
  $Q_{1}$   &$73.13$            &$Q_{11}$   &$33310.28$        &$Q_{21}$   &$228232.34$   \\
  $Q_{2}$   &$319.17$           &$Q_{12}$   &$43100.48$        &$Q_{22}$   &$262277.55$  \\
  $Q_{3}$   &$851.80$           &$Q_{13}$   &$54654.69$        &$Q_{23}$   &$299556.76$ \\
  $Q_{4}$   &$1822.69$          &$Q_{14}$   &$68119.90$        &$Q_{24}$   &$340216.96$\\
  $Q_{5}$   &$3381.01$          &$Q_{15}$   &$83643.10$        &$Q_{25}$   &$384405.17$ \\
  $Q_{6}$   &$5674.24$          &$Q_{16}$   &$101371.31$       &$Q_{26}$   &$432268.38$ \\
  $Q_{7}$   &$8849.45$          &$Q_{17}$   &$121451.52$       &$Q_{27}$   &$483953.58$ \\
  $Q_{8}$   &$13053.65$         &$Q_{18}$   &$144030.72$       &$Q_{28}$   &$539607.79$ \\
  $Q_{9}$   &$18433.86$         &$Q_{19}$   &$169255.93$       &$Q_{29}$   &$599378.00$ \\
  $Q_{10}$  &$25137.07$         &$Q_{20}$   &$197274.14$       &$Q_{30}$   &$663411.21$ \\
  \hline
\end{tabular}
\end{table}

\begin{thm} Suppose $Q_n$ is a fractal M\"{o}bius octagonal networks of length $n\geq2$. Therefore,
\begin{eqnarray*}
Kc(Q_n)&=&\sum_{j=2}^{j=3n}\frac{1}{\alpha_j}+\sum_{j=1}^{j=3n}\frac{1}{\rho_j}\\
&=&\frac{147n^{2}-19}{84}+\frac{\frac{37\sqrt{15}n}{30}\left[(\frac{4+\sqrt{15}}{12})^{n}-(\frac{4-\sqrt{15}}{12})^{n}\right]}
{(\frac{4+\sqrt{15}}{12})^{n}+(\frac{4-\sqrt{15}}{12})^{n}+2(\frac{1}{12})^{n}}.
\end{eqnarray*}
\end{thm}

\begin{thm} Suppose $Q_n$ is a fractal M\"{o}bius octagonal networks of length $n\geq2$. Then
\begin{eqnarray*}
\tau(Q_n)&=&\frac{3n}{2}[(4+\sqrt{15})^{n}+(4-\sqrt{15})^{n}+2].
\end{eqnarray*}
\end{thm}
\noindent\textbf{Proof.} According to Fact 4, one has
$$\prod_{j=2}^{3n}\alpha_j=(-1)^{3n-1}d_{3n-1}=21n^{2}(\frac{1}{12})^{n}.$$\\
Similarly, by Claim 4, one finds\\
\begin{eqnarray*}
\prod_{j=1}^{3n}\rho_j=det  \mathcal
L_S=(\frac{4+\sqrt{15}}{12})^{n}+(\frac{4-\sqrt{15}}{12})^{n}+2(\frac{1}{12})^{n}.
\end{eqnarray*}
It needs to be pointed out
$\prod_{j=1}^{6n}d_j(Q_n)=2^{4n}3^{2n},~and~|E(Q_n)|=7n.$ Then Theorem 3.5 is acquired by the formula of
$\tau(Q_n)$.

Finally, according to the formula, we calculate the number of spanning trees from $Q_1$ to $Q_{12}$.
%\noindent{\bf Proof.} Based on the Lemma 4.1, one immediately gets
%\begin{eqnarray*}
%Kf^*(O_n)&=&2\cdot(13n+1)\cdot(\sum_{i=2}^{3n+1}\frac{1}{\lambda_i}+\sum_{j=1}^{3n+1}\frac{1}{\delta_j})\\
%&=&2(13n+1)\cdot\bigg(\frac{169n^3+39n^2+22n}{4(13n+1)}+\frac{17n+4}{6}\bigg)\\
%&=&\frac{507n^3+559n^2+204n+8}{6}.
%\end{eqnarray*}
%The proof of the theorem completed.\hfill\rule{1ex}{1ex}
\begin{table}[h]
\setlength{\abovecaptionskip}{0.05cm} \centering\vspace{.3cm}
\caption{The complexity of $Q_n$ from $Q_{1}$
to $Q_{12}$.}
\begin{tabular}{cccccccc}
  \hline
  $G$ &  ~~ $\tau(G)$ ~~  & $G$ & ~~$\tau(G)$   ~~& $G$ & ~~$\tau(G)$ & ~~$G$ & ~~$\tau(G)$ \\
  \hline
  $Q_{1}$ & $15$ & $Q_{4}$ & $230,64$ & ~~$Q_{7}$ & $196,863,45$ & ~~$Q_{10}$ & $137,241,225,60$  \\
  $Q_{2}$ & $192$ & $Q_{5}$ & $226,875$ & ~~$Q_{8}$ & $177,131,568$ & ~~$Q_{11}$ & $118,854,766,965$
  \\
  $Q_{3}$ & $2205$ & $Q_{6}$ & $214,329,6$ & ~~$Q_{9}$ & $156,887,293,5$ &~~ $Q_{12}$ & $102,080,901,875,2$
  \\
  \hline
\end{tabular}
\end{table}

\section{Conclusion} \ \ \ \ \
Under the research of some scholars, we have carried on some expansion, studied a new graph, fractal
M\"{o}bius oOctagonal networks $(Q_n)$. In this paper, we first restate the normalized Laplacian
decomposition theorem. Then, the product of the sum of reciprocal eigenvalues of $\mathcal L_A$ and
$\mathcal L_S$ are required, by using the Vieta's theorem for the characteristic polynomials of $\mathcal
L_A$ and $\mathcal L_S$. Finally, $DK(Q_n),~Kc(Q_n)$ and $\tau(Q_n)$ of fractal M\"{o}bius octagonal
networks $(Q_n)$ are obtained.
%\ \ \ \ The authors would like to express their
%sincere gratitude to the editor and referees for a very careful
%reading of the paper and for all their insightful comments and
%valuable suggestions, which led to a number of improvements in
%this paper
%, Anhui
%Provincial Natural Science Foundation (Nos. KJ2015A331,
%KJ2013B105)
%, and the Special Fund for Basic Scientific Research of
%Central Colleges, South-Central University for
%Nationalities(CZY18032), and China Postdoctoral Science Foundation
%under Grant no. 2017M621579, the Postdoctoral Science Foundation
%of Jiangsu Province under Grant no. 1701081B, Project of Anhui
%Jianzhu University under Grant no. 2016QD116 and 2017dc03, Anhui
%Province Key Laboratory of Intelligent Building \& Building Energy
%Saving.

%\section*{Conflicts of Interest} The authors declare that there are no conflicts of interest regarding the
publication of this paper.

%\section*{ORCID}
%\Author{Jia-Bao Liu $^{1}$\orcidA{}}.
\section*{Funding}

This work was funded in part by Anhui Provincial Natural Science Foundation under Grant 2008085J01, and by
National Natural Science Foundation of China Grant 11601006,
and by China Post-doctoral Science Foundation under Grant 2017M621579.


\begin{thebibliography}{99}
\small \setlength{\itemsep}{-.8mm}
\bibitem{D.J} D.J. Watts, S.H. Strogatz, Collective dynamics of `small-world' networks, Nature. \textbf{393}
    (1998) 440-442.

\bibitem{A.L} A.-L. Barab\'{a}si, R. Albert, Emergence of scaling in random networks, Science. \textbf{286}
    (1999) 509-512.

\bibitem{WA} S. Wang, L. Xi , H. Xu, Li. Wang, Scale-free and small-world properties of Sierpinski networks,
    Physica A.  \textbf{465} (2017) 690-700.

\bibitem{J.A} J. A. Bondy, U. S. R. Murty, Graph theory, Springer, New York, 2008.

\bibitem{b7} D.J. Klein, M. Randi\'{c}, Resistance distance, J. Math. Chem. \textbf{12}(1) (1993) 81-95.

\bibitem{b8} D.J. Klein, Resistance-distance sum rules, Croat. Chem. Acta. \textbf{75}(2) (2002) 633-649.

\bibitem{b9} D.J. Klein, O. Ivanciuc, Graph cyclicity, excess conductance, and resistance deficit, J. Math.
    Chem. \textbf{30}(3) (2001) 271-287.

\bibitem{b10} I. Gutman, B. Mohar, The quasi-Wiener and the Kirchhoff indices coincide, J. Chem. Inf.
    Comput. Sci. \textbf{36}(5) (1996) 982-985.

\bibitem{b11}  H. Y. Zhu, D. J. Klein,  I. Lukovits, Extensions of the Wiener Number, J. Chem. Inf. Comput.
    Sci. \textbf{36}(3) (1996) 420-428.

\bibitem{b34} H. Y. Chen, F. J. Zhang, Resistance distance and the normalized Laplacian spectrum, Discrete
    Appl. Math. \textbf{155}(5) (2007) 654-661.

\bibitem{b35} L. H. Feng, I. Gutman, G .H. Yu, Degree Kirchhoff index of unicyclic graphs, MATCH Commun.
    Math. Comput. Chem. \textbf{69} (2013) 629-648.

\bibitem{b36} J. Huang, S. C. Li, On the normalised Laplacian spectrum, degree-Kirchhoff index and spanning
    trees of graphs, Bull. Aust. Math. Soc. \textbf{91} (2015) 353-367.

\bibitem{b13} S. Butler, A note about cospectral graphs for the adjacency and normalized Laplacian matrices,
    Linear Multilinear Algebra \textbf{58}(3) (2010) 387-390.

\bibitem{b14} M. S. Cavers, S. Fallat, S. Kirkland, On the normalized Laplacian energy and general
    Randi\'{c} index $R_{-1}$ of graphs, Linear Algebra Appl. \textbf{433}(1) (2010) 172-190.

\bibitem{b16}  H. Chen, J. Jost, Minimum vertex covers and the spectrum of the normalized Laplacian on
    trees, Linear Algebra Appl. \textbf{437}(4) (2012) 1089-1101.

\bibitem{b17} K. Ch. Das, S. W. Sun, Normalized Laplacian eigenvalues and energy of trees, Taiwanese J.
    Math. \textbf{20}(3) (2016) 491-507.

\bibitem{b15}  G.T. Chen, D. George, H. Frank, Z.S. Li, P. Kinnari, S. Michael, An interlacing result on
    normalized Laplacians, SIAM J. Discrete Math. \textbf{18}(2) (2004) 353-361.

\bibitem{b18} M. Levene, G. Loizou, Kemeny's costant and the random surfer, Amer. Math. Monthly.
    \textbf{109}(8) (2002) 741-745.

\bibitem{b25}J. Huang, S. C. Li, L. Q. Sun, The normalized Laplacians, degree-Kirchhoff index and the
    spanning trees of linear hexagonal chains, Discrete Appl. Math.
\textbf{207}(10) (2016) 67-79.

\bibitem{b19} Q. Li,  S. Zaman, W. Sun, J.  Alam, Study on the normalized Laplacian of a penta-graphene with
    applications, Int. J. Quantum Chem. \textbf{120}(9) (2020) e26154.

\bibitem{b20} S. Li, W. Wei, S. Yu, On normalized Laplacians, multiplicative degree-Kirchhoff indices, and
    spanning trees of the linear [n] phenylenes
and their dicyclobutadieno derivatives, Int. J. Quantum Chem. \textbf{119} (2019) e25863.

\bibitem{b21} X. Ma, H. Bian, The normalized Laplacians, degree-Kirchhoff index and the spanning trees of
    cylinder phenylene chain, Polycycl. Aromat. Comp. 2019. DOI: 10.1080/10406638.2019.1665553

\bibitem{b22} X. Ma, H. Bian, The normalized Laplacians, degree-Kirchhoff index and the spanning trees of
    hexagonal M\"{o}bius graphs, Appl. Math. Comput. \textbf{355} (2019) 33-46.



\bibitem{b23} J. B. Liu, J. Zhao, Z. X. Zhu, On the number of spanning trees and normalized Laplacian of
    linear octagonal-quadrilateral networks,
Int. J. Quantum Chem. \textbf{119} (2019) e25971.

\bibitem{b24} Q. Zhu, Kirchhoff index, degree-Kirchhoff index and spanning trees of linear octagonal chains,
    Australas. J. Comb. \textbf{118} (2018) e25787.

\bibitem{b26} Z. Zhu, J. B. Liu, The normalized Laplacian, degree-Kirchhoff index and the spanning tree
    numbers of generalized phenylenes, Discrete Appl. Math. \textbf{254}(15) (2019) 256-267.

\bibitem{b27} C. Liu, Y. Pan, J. Li, On the Laplacian spectrum and Kirchhoff index of generalized
    phenylenes, Polycycl. Aromat. Comp. 2019. DOI: 10.1080/10406638.2019.1703765.

\bibitem{b28} C. He, S. Li, W. Luo, Calculating the normalized Laplacian spectrum and the number of spanning
    trees of linear pentagonal chains, J. Comput. Appl. Math. \textbf{344}(15) (2018) 381-393.



\bibitem{b30} H. Y. Chen, F. J. Zhang, Resistance distance and the normalized Laplacian spectrum, Discrete
    Appl. Math. \textbf{155}(5) (2007) 654-661.

\bibitem{b32} S. Butler, Algebraic aspects of the normalized Laplacian, in: A. Beveridge, J. Griggs, L.
    Hogben, G. Musiker, e. P. Tetali (Eds.), Recent Trends in Combinatorics, The IMA Volumes in Mathematics
    and its Applications, IMA. \textbf{159} (2016) 295-315.

\bibitem{b12} F. R. K. Chung, Spectral Graph Theory, American Mathematical Society, Providence, RI 1997.

%\bibitem{J.-B} J.-B. Liu, X.F. Pan, L. Yu, D. Li, Complete characterization of bicyclic graphs with minimal
Kirchhoff index, Discret. Appl. Math. 200 (2016)
%95-107.
%
%\bibitem{J.-B.} J.-B. Liu, X.F. Pan, Minimizing Kirchhoff index among graphs with a given vertex
%bipartiteness, Appl. Math. Comput. 291 (2016) 84-88.
%
%
%
%
%
%\bibitem{A.D} A. Dobrynin, Branchings in trees and the calculation of the Wiener index of a tree, MATCH
Commun. Math. Comput. Chem. 41 (2000) 119-134.
%
%
%
%
%
%
%
%
%
%
%
%
%
%\bibitem{H.Y} H. Y. Zhu, D. J. Klein, I. Lukovits, Extensions of the Wiener number, J. Chem. Inf. Comput.
Sci. 36 (1996) 420-428.


\end{thebibliography}
\end{document}